\RequirePackage{fix-cm}

\documentclass[smallextended]{svjour3}       

\smartqed

\usepackage{graphicx}
\usepackage{color}
\usepackage{amsmath}
\usepackage{amssymb}
\usepackage{amsfonts}
\usepackage{subcaption}
\usepackage{xcolor}
\usepackage{tabularx,booktabs}
\usepackage{listings}
\usepackage{tikz}
\usepackage{pgfkeys}
\usepackage{adjustbox}
\usepackage{footnotehyper}
\usepackage{placeins}
\usetikzlibrary{shapes,arrows}
\usepackage{enumitem}
\usepackage{multirow}
\usepackage[colorlinks=true]{hyperref}
\hypersetup{urlcolor=blue, citecolor=blue}
\usepackage{breakurl}
\usepackage{pifont}

\lstdefinestyle{style1}{
    basicstyle=\ttfamily,
    keywordstyle =    \color{black}, 
    keywordstyle =    \color{black}, 
    commentstyle =    \color{black},
    stringstyle  =    \color{black},
    columns=fullflexible,
    keepspaces=true,
    upquote=true,
}

\usepackage{pgfplots}
\pgfplotsset{compat=1.12}

\tikzstyle{loosely dashed}=          [dash pattern=on 2pt off 4pt]

\pgfplotsset{plotOptions4/.style={%
		width=\linewidth,
				ymax=1.1,ymin=.0001,
				xmin=0.9,xmax=100,
		xlabel={Iteration count (log scale)},
		ylabel={Worst-case guarantee},
		label style={font=\footnotesize},
		legend style={font=\scriptsize},
		ytick={1,.1, 0.01, 0.001, 0.0001},
		xtick={1, 10, 50},
		tick label style={font=\footnotesize},
		solid,
		very thick
	}}
	
\pgfplotsset{plotOptions3/.style={%
		width=\linewidth,
				ymax=1.05,ymin=.85,
				xmin=10,xmax=100000,
		xlabel={Condition number},
		ylabel={Worst-case guarantee},
		label style={font=\footnotesize},
		legend style={font=\scriptsize},
		ytick={1.05,1,0.95,0.9,0.85},
		xtick={10,100,1000,10000,100000,1000000},
		tick label style={font=\footnotesize},
		solid,
		very thick
	}}

\pgfplotsset{plotOptions2/.style={%
		width=\linewidth,
				ymax=10,ymin=0.001,
				xmin=1,xmax=100,
		xlabel={Iteration count (log scale)},
		ylabel={Worst-case guarantee},
		label style={font=\footnotesize},
		legend style={font=\scriptsize},
		ytick={10, 1, 0.1, 0.1, 0.01, 0.001},
		xtick={1, 10, 100},
		tick label style={font=\footnotesize},
		solid,
		very thick
	}}
	
\pgfplotsset{plotOptions1/.style={%
		width=\linewidth,
				ymax=1,ymin=0.0001,
				xmin=0.9,xmax=100,
		xlabel={Iteration count (log scale)},
		ylabel={Worst-case guarantee},
		label style={font=\footnotesize},
		legend style={font=\footnotesize},
		ytick={1, 0.1, 0.1, 0.01, 0.001, 0.0001},
		xtick={1, 10, 50},
		tick label style={font=\footnotesize},
		solid,
		very thick
	}}

\pgfplotsset{plotOptions5/.style={%
		width=\linewidth,
				ymax=1,ymin=0.,
				xmin=0,xmax=2,
		xlabel={Step-size},
		ylabel={Worst-case guarantee},
		label style={font=\footnotesize},
		legend style={font=\footnotesize},
		ytick={1, 0.25, 0.5, 0.75, 1},
		xtick={0, 0.5, 1, 1.5, 2},
		tick label style={font=\footnotesize},
		solid,
		very thick
	}}

\newcommand{\Fmul}{\mathcal{F}_{\mu,L}}
\newcommand{\Rd}{\mathbb{R}^d}

\newcommand{\tr}[1]{\mathrm{Tr}\left(#1\right)}
\renewcommand{\leq}{\leqslant}
\renewcommand{\geq}{\geqslant}
\renewcommand{\succeq}{\succcurlyeq}

\newcommand{\eqdef}{\triangleq}

\journalname{Mathematical Programming Computation}

\DeclareMathOperator*{\argmin}{arg\,min}
\newcommand{\Span}{\mathrm{span}}
\usepackage{stmaryrd}
\newcommand{\range}[2]{\llbracket #1, #2 \rrbracket}

\begin{document}

\title{\textsc{PEPit}: computer-assisted worst-case analyses of first-order optimization methods in Python}
\subtitle{}

\titlerunning{\textsc{PEPit}: computer-assisted worst-case analyses in Python}

\author{Baptiste Goujaud
        \and
        C\'eline Moucer
        \and
        Fran\c{c}ois Glineur
        \and
        Julien M. Hendrickx
        \and
        Adrien B. Taylor
        \and
        Aymeric Dieuleveut
}

\institute{
Baptiste Goujaud \at CMAP, École Polytechnique, Institut Polytechnique de Paris, Paris, France. \\
\email{baptiste.goujaud@polytechnique.edu}
\and
C\'eline Moucer \at INRIA \& D.I. \'Ecole Normale Supérieure,
CNRS \& PSL Research University, Paris, France, Ecole des Ponts, Marne-la-Vallée, France. \\
\email{celine.moucer@inria.fr}
\and
Fran\c{c}ois Glineur \at ICTEAM institute \& CORE, UCLouvain, Louvain-la-Neuve, Belgium. \\
\email{francois.glineur@uclouvain.be}
\and
Julien M. Hendrickx\at ICTEAM institute, UCLouvain, Louvain-la-Neuve, Belgium. \\
\email{julien.hendrickx@uclouvain.be}
\and
Adrien B. Taylor \at INRIA \& D.I. \'Ecole Normale Supérieure,
CNRS \& PSL Research University, Paris, France \\
\email{adrien.taylor@inria.fr}
\and
Aymeric Dieuleveut \at CMAP, École Polytechnique, Institut Polytechnique de Paris, Paris, France. \\
\email{aymeric.dieuleveut@polytechnique.edu}
}

\date{Received: date / Accepted: date}

\maketitle

\begin{abstract}
    \textsc{PEPit} is a \textsc{python} package aiming at simplifying the access to worst-case analyses of a large family of first-order optimization methods possibly involving gradient, projection, proximal, or linear optimization oracles, along with their approximate, or Bregman variants.

    In short, \textsc{PEPit} is a package enabling computer-assisted worst-case analyses of first-order optimization methods. The key underlying idea is to cast the problem of performing a worst-case analysis, often referred to as a performance estimation problem (PEP), as a semidefinite program (SDP) which can be solved numerically. To do that, the package users are only required to write first-order methods nearly as they would have implemented them. The package then takes care of the SDP modeling parts, and the worst-case analysis is performed numerically via a standard solver.
\end{abstract}

\section{Introduction}

Due to their low cost per iteration, first-order optimization methods became a major tool in the modern numerical optimization toolkit. Those methods are particularly well suited when targeting only low to medium accuracy solutions, and play a central role in many fields of applications that include machine learning and signal processing. Their simplicity further allows both occasional and expert users to use them. On the contrary, when it comes to their analyses (usually based on worst-case scenarios), they are mostly reserved for expert users. The main goal of this work is to allow simpler and reproducible access to worst-case analyses for first-order methods.

\textsc{PEPit} is a \textsc{python} package enabling computer-assisted worst-case analysis of a large family of first-order optimization methods. After being provided with a first-order method and a standard problem class, the package reformulates the problem of performing a worst-case analysis as a semidefinite program~(SDP). This technique is commonly referred to as \emph{performance estimation problems}~(PEPs) and was introduced by~\cite{drori2014performance,drori2014contributions}. The package uses PEPs as formalized by~\cite{taylor2017smooth,taylor2017exact}.

In short, performing a worst-case analysis of a first-order algorithm usually relies on four main ingredients: a first-order algorithm (to be analyzed), a class of problems (containing the assumptions on the function to be minimized), a performance measure (measuring the quality of the output of the algorithm under consideration; for convenience here we assume that the algorithm aims at minimizing this performance measure and our analysis aims at finding a worst-case guarantee on it), and an initial condition (measuring the quality of the initial iterate). Performing the worst-case analysis (i.e., computing worst-case scenarios) corresponds to maximizing the performance measure of the algorithm on the class of problems, under a constraint induced by the initial condition. It turns out that such optimization problems can often be solved using SDPs in the context of first-order methods.

PEPs provide a principled approach to worst-case analyses but usually rely on potentially tedious semidefinite programming (SDP) modeling and coding steps. 
The \textsc{PEPit} package eases access to the methodology by automatically handling the modeling part, thereby limiting the amount of time spent on this tedious task and the risk of introducing coding mistakes in the process. In short, this work allows users to (i)~write their first-order algorithms nearly as they would have implemented them, (ii)~let \textsc{PEPit} (a) perform the modeling and coding steps, and (b) perform the worst-case analysis numerically using tools for semidefinite programming in \textsc{python}~\cite{mosek,diamond2016cvxpy,o2016conic}.

As a result, the package enables users to easily obtain worst-case analyses for most of the standard first-order methods, classes of problems, performance measures, and initial conditions. 
This is useful to numerically verify existing convergence guarantees, as well as to ease the development of new analyses and methods.
To this end, the toolbox contains tools for analyzing classical scenarios of the first-order literature: standard problem classes (such as convex functions, smooth convex functions, Lipschitz convex functions, etc.) and algorithmic operations (such as gradient, proximal, or linear optimization oracles, etc.). Finally, the package contains more than $75$ examples and is designed in an open fashion, allowing users to easily add new ingredients (such as their own problem classes, oracles, or algorithms as examples).

\paragraph{Organization of the paper.} \label{p:content} This paper is organized as follows.
First, Section~\ref{s:example} exemplifies the PEP approach on a very simple example, namely computing a worst-case contraction factor for gradient descent, and shows how to code this example in \textsc{PEPit}.
Section~\ref{s:SDP_details} provides details on the semidefinite programs that can be formulated through the package along with the relationship between those formulations and the coding steps.
Then, Section~\ref{s:content} provides a roadmap through the package.
Finally, Section~\ref{sec:numerical-examples} provides three additional numerical examples (including a composite minimization problem and a stochastic one), and some concluding remarks and perspectives are drawn in Section~\ref{s:ccl}.

\paragraph{Related works.} The \textsc{PEPit} package relies on performance estimation problems as formalized in~\cite{taylor2017exact}.
It also contains some improvements and generalizations to other problem and algorithmic classes such as monotone and nonexpansive operators~\cite{ryu2020operator,lieder2020convergence}, quadratic optimization~\cite{bousselmi2023interpolation}, stochastic methods and verification of potential (or Lyapunov/energy) functions~\cite{hu2018dissipativity,fazlyab2018analysis,taylor19bach} as inspired by the related control-theoretic IQC framework~\cite{lessard2016analysis}.
The package also contains numerous examples; e.g., recent analyses and developments from~\cite{kim2016optimized,van2017fastest,gu2020tight,kim2018optimizing,kim2019accelerated,lieder2020convergence,gannot2021frequency,drori2021optimal,abbaszadehpeivasti2021exact,gorbunov2022extragradient}.
The package can be seen as an extended open source \textsc{python} version of the \textsc{matlab} package \textsc{pesto}~\cite{taylor2017performance} on various aspects such as its documentation, its coding style and its access through standard open interfaces (such as pip), \textsc{PEPit} is more professional than \textsc{pesto}.

\paragraph{Dependencies{.}} The package heavily builds on existing software for solving semidefinite programs, including CVXPY~\cite{diamond2016cvxpy}, SCS~\cite{o2016conic}, and MOSEK~\cite{mosek}.

\section{\textsc{PEPit} on a simple example}\label{s:example}
In this section, we illustrate the use of the package for studying the worst-case properties of a standard scenario: gradient descent for minimizing a smooth strongly convex function.
The goal of this elementary example is twofold.
First, we want to provide the base mathematical steps enabling the use of semidefinite programming for performing worst-case analyses, together with a corresponding \textsc{PEPit} code.
Second, we want to highlight the main ingredients that can be generalized to other problem setups (e.g., Theorem~\ref{thm:interp} below providing ``interpolation conditions'' for the class of smooth strongly convex functions), allowing us to analyze more algorithms under different assumptions (which are listed in Section~\ref{s:content}).

For this example, we consider the convex optimization problem
\begin{equation}
    \label{eq:min_problem} \min_{x\in\mathbb{R}^d} f(x),
\end{equation}
where $f$ is $L$-smooth and $\mu$-strongly convex (notation $f\in\Fmul(\mathbb{R}^d)$, or $f\in\Fmul$ when $d$ is unspecified). So we assume $f$ to satisfy
\begin{itemize}
    \item ($L$-smoothness) $\forall x,y\in\mathbb{R}^d$ we have that 
    \[ f(x) \leqslant f(y)+\langle \nabla f(y);x-y\rangle+\frac{L}{2}\|x-y\|^2_2,\]
    \item ($\mu$-strong convexity) $\forall x,y\in\mathbb{R}^d$ we have that 
    \[ f(x) \geqslant f(y)+\langle \nabla f(y);x-y\rangle+\frac{\mu}{2}\|x-y\|^2_2.\]
\end{itemize}
Our goal for the rest of this section is to show how to compute the smallest possible~$\tau(\mu, L, \gamma)$ (often referred to as the ``contraction factor'') such that
\begin{equation}
     \| x_{1}-y_{1}\|_2^2\leqslant \tau(\mu,L,\gamma) \|x_0-y_0\|^2_2,\label{eq:tau_GD}
\end{equation}
is valid for all $f\in\Fmul$ and all $x_0,y_0\in\mathbb{R}^d$ when $x_{1}$ and $y_{1}$ are obtained from gradient steps from respectively $x_0$ and $y_0$. That is, $x_{1}=x_0-\gamma \nabla f(x_0)$ and $y_{1}=y_0-\gamma \nabla f(y_0)$. First, we show that the problem of computing $\tau(\mu, L, \gamma)$ can be framed as a semidefinite program~(SDP), and then illustrate how to use \textsc{PEPit} for computing it without going into the SDP modeling details.

\subsection{A performance estimation problem for the gradient method}\label{s:pep}
It is relatively straightforward to establish that the smallest possible $\tau(\mu,L,\gamma)$ for which~\eqref{eq:tau_GD} is valid can be computed as the worst-case value of $\| x_{1}-y_{1}\|_2^2$ when $\|x_0-y_0\|^2_2\leqslant 1$. That is, we compute $\tau(\mu,L,\gamma)$ as the optimal value to the following optimization problem:
\begin{equation}\label{eq:PEP1}
    \begin{aligned}
    \tau(\mu,L,\gamma)=\max_{\substack{f,d\\x_0,x_{1}\in\Rd\\y_0,y_{1}\in\Rd}} &\| x_{1}-y_{1}\|_2^2 \\
    \text{s.t.}\,\, & d\in\mathbb{N},\, f\in\Fmul(\Rd),\\
    & \| x_{0}-y_{0}\|_2^2\leqslant 1,\\
    & x_{1}=x_0-\gamma \nabla f(x_0),\\
    & y_{1}=y_0-\gamma \nabla f(y_0).
    \end{aligned}
\end{equation}
As written in~\eqref{eq:PEP1}, this problem involves an infinite-dimensional variable $f$. Our first step towards formulating~\eqref{eq:PEP1} as an SDP consists of reformulating it by sampling $f$ (i.e., evaluating its function value and gradient) at the two points where its gradient is evaluated:
\begin{equation}\label{eq:PEP_sampled}
    \begin{aligned}
    \tau(\mu,L,\gamma)=\max_{\substack{d,f_{x_0},f_{y_0}\\x_0,x_{1},g_{x_0}\in\Rd\\y_0,y_{1},g_{y_0}\in\Rd}} &\| x_{1}-y_{1}\|_2^2 \\
    \text{s.t.}\,\, & d\in\mathbb{N},\\
    & \| x_{0}-y_{0}\|_2^2\leqslant 1,\\
    &\exists f\in\Fmul(\Rd): \left\{\begin{array}{cc}
    f(x_0)=f_{x_0} & \nabla f(x_0)=g_{x_0}  \\
    f(y_0)=f_{y_0} & \nabla f(y_0)=g_{y_0}
    \end{array}\right.\\
    & x_{1}=x_0-\gamma g_{x_0},\\
    & y_{1}=y_0-\gamma g_{y_0},
    \end{aligned}
\end{equation}
where we replaced the variable $f$ by its discrete version, which we constrain to be ``interpolable'' (or ``extendable'') by a smooth strongly convex function over $\mathbb{R}^d$. To arrive at a tractable problem, we use the following interpolation (or extension) result.
\begin{theorem}\cite[Theorem 4]{taylor2017smooth}\label{thm:interp} Let $I$ be an index set and $S=\{(x_i,g_i,f_i)\}_{i\in I}$ be such that $x_i,g_i\in\mathbb{R}^d$ and $f_i\in\mathbb{R}$ for all $i\in I$. There exists a function $F\in\mathcal{F}_{\mu,L}(\mathbb{R}^d)$ such that $f_i=F(x_i)$ and $g_i=\nabla F(x_i)$ (for all $ i\in I$) if and only if for all $i,j\in I$ we have
\begin{align}
f_i\geq f_j+\langle g_j;x_i-x_j\rangle+\tfrac{1}{2L}\|g_j-g_i\|^2_2+\tfrac{\mu L}{2(L-\mu)}\|x_i-x_j-\tfrac1L(g_i-g_j)\|^2_2.\label{eq:muLinterp}
\end{align}
\end{theorem}\vspace{.5cm}
Using Theorem~\ref{thm:interp}, we can formulate the problem of computing $\tau(\mu,L,\gamma)$ as a (nonconvex) quadratic problem:
\begin{equation}\label{eq:PEP_nncvx_QP}
    \begin{aligned}
    \max_{\substack{d,f_{x_0},f_{y_0}\\x_0,g_{x_0}\in\Rd\\y_0,g_{y_0}\in\Rd}} &\| (x_0-\gamma g_{x_0})-(y_0-\gamma g_{y_0})\|_2^2 \\
    \text{s.t.}\,\, & d\in\mathbb{N},\\
    & \|x_{0}-y_{0}\|_2^2\leqslant 1,\\
    & f_{y_0}\geq f_{x_0}+\langle g_{x_0};y_0-x_0\rangle+\tfrac{1}{2L}\|g_{y_0}-g_{x_0}\|^2_2 \\ 
    & \hspace{1.42cm} +\tfrac{\mu L}{2(L-\mu)}\|x_0-y_0-\tfrac1L(g_{x_0}-g_{y_0})\|^2_2,\\
    & f_{x_0}\geq f_{y_0}+\langle g_{y_0};x_0-y_0\rangle+\tfrac{1}{2L}\|g_{y_0}-g_{x_0}\|^2_2 \\
    & \hspace{1.42cm} +\tfrac{\mu L}{2(L-\mu)}\|x_0-y_0-\tfrac1L(g_{x_0}-g_{y_0})\|^2_2.
    \end{aligned}
\end{equation}
Relying on a standard trick from semidefinite programming, one can convexify this problem using a Gram representation of the variable (this is due to maximization over $d$). That is, we formulate the problem using a positive semidefinite matrix $G\succcurlyeq 0$ defined as
\begin{equation}\label{eq:def_G}
    G\eqdef \begin{pmatrix} \|x_0-y_0\|^2_2 & \langle x_0-y_0; g_{x_0}\rangle & \langle x_0-y_0; g_{y_0}\rangle\\
\langle x_0-y_0; g_{x_0}\rangle & \|g_{x_0}\|_2^2 & \langle g_{x_0}; g_{y_0}\rangle\\
\langle x_0-y_0; g_{y_0}\rangle & \langle g_{x_0}; g_{y_0}\rangle & \|g_{y_0}\|_2^2
\end{pmatrix}\succcurlyeq 0.
\end{equation}
Using this change of variable, we arrive {at}
\begin{equation}\label{eq:PEP_SDP}
    \begin{aligned}
    \max_{f_{x_0},f_{y_0},G} &G_{1,1}-2\gamma (G_{1,2}-G_{1,3})+ \gamma^2 (G_{2,2}+G_{3,3}-2G_{2,3}) \\
    \text{s.t.}\,\, & G\succcurlyeq 0,\\
    &G_{1,1}\leqslant 1,\\
    &f_{y_0}\geq f_{x_0}+\tfrac{1}{L-\mu}\left(\tfrac{\mu L}{2} G_{1,1}-LG_{1,2}+\mu G_{1,3}+\tfrac{1}{2}G_{2,2}-G_{2,3}+\tfrac{1}{2}G_{3,3}\right),\\
    &f_{x_0}\geq f_{y_0}+\tfrac{1}{L-\mu}\left(\tfrac{\mu L}{2}G_{1,1}-\mu G_{1,2}+L G_{1,3}+\tfrac{1}{2}G_{2,2}-G_{2,3}+\tfrac{1}{2}G_{3,3}\right),
    \end{aligned}
\end{equation}
which can be solved numerically using standard tools, see, e.g.,~\cite{diamond2016cvxpy,mosek,o2016conic}. Using numerical and/or symbolical computations, one can then easily arrive at $\tau(\mu, L,\gamma)=\max\{(1-L\gamma)^2,(1-\mu \gamma)^2\}$ and hence that
\begin{equation}
     \| x_{1}-y_{1}\|_2^2\leqslant  \max\{(1-L\gamma )^2,(1-\mu \gamma)^2\}\|x_0-y_0\|^2_2,\label{eq:tau_GD_ex}
\end{equation}
for all $d\in\mathbb{N}$, $f\in\Fmul(\Rd)$ and $x_0,y_0\in\mathbb{R}^d$ when $x_{1},y_{1}\in\mathbb{R}^d$ are generated from gradient steps from respectively $x_0$ and $y_0$. In the next section, we show how to perform this analysis using \textsc{PEPit}, which automates the sampling {(i.e., the evaluations of the function or its gradients on given points)} and SDP-modeling procedures. In more complex settings where more functions need to be sampled and/or more iterates have to be taken into account, avoiding those steps allows to largely limits the probability of making a mistake in the process of performing the worst-case analysis (numerically) while allowing to spare a significant amount of time in the process.

\begin{remark}[Important ingredients for the SDP reformulations]
To understand what \textsc{PEPit} can do, it is crucial to understand which elements allowed to cast the worst-case analysis as such a semidefinite program (which is what we refer to as the ``modeling'' of the problem). In short, the SDP reformulation of the worst-case computation problem was made possible due to $4$ main ingredients (see, e.g.,~\cite[Section 2.2]{taylor2017exact}):
\begin{enumerate}[noitemsep]
    \item the algorithmic steps can be expressed linearly in terms of the iterates and gradient values (i.e., step-sizes do not depend on the function at hand),
    \item the class of functions has ``interpolation condition'' \footnote{Interpolation conditions characterize the \textbf{existence} of a function (that has particular function values and gradients at given points) in the considered class by a list of constraints on those gradients, points, and function values.
    Such interpolation theorems (see, e.g., Theorem~\ref{thm:interp}) have been obtained in the literature for various problem classes, see, e.g.,~\cite{taylor2017smooth,taylor2017exact,ryu2020operator}.} that are linear in $G$ and the function values, 
    \item the performance measure is linear (or convex piecewise linear) in $G$ and the function values,
    \item the initial condition is linear in $G$ and the function values.
\end{enumerate}
Those ingredients allow the use of PEPs much beyond the simple setup of this section.
That is, PEPs apply for performing worst-case analyses involving a variety of first-order oracles, initial conditions, performance measures, and problem classes (see Section~\ref{s:SDP_details} for the general modeling of the problem, and Section~\ref{s:content} for a non-exhaustive list of cases that are covered).
\end{remark}

\subsection{Code}\label{s:code}

In the previous section, we introduced the PEP and SDP modeling steps for computing a tight worst-case contraction factor for gradient descent in the form~\eqref{eq:tau_GD}. Although this particular SDP~\eqref{eq:PEP_SDP} might be solved analytically, many optimization methods lead to larger SDPs with more complicated structures.
In general, we can reasonably only hope to solve them numerically. In the following lines, we describe how to use \textsc{PEPit} for computing a \emph{contraction factor} without explicitly going into the modeling steps. Compared to previous section, we allow ourselves to perform $n\in\mathbb{N}$ iterations and compute the smallest possible value of $\tau(\mu,L,\gamma,n)$ such that
\begin{equation}
     \| x_{n}-y_{n}\|_2^2\leqslant \tau(\mu,L,\gamma,n) \|x_0-y_0\|^2_2,\label{eq:tau_GD_n}
\end{equation}
where $x_n$ and $y_n$ are computed from $n$ iterations of gradient descent with step-size $\gamma$ starting from respectively $x_0$ and $y_0$. As illustrated in the previous section for the case $n=1$, computing the smallest possible such $\tau(\mu,L,\gamma,n)$ is equivalent  to computing the worst-case value of $\| x_{n}-y_{n}\|_2^2$ under the constraint that $\| x_{0}-y_{0}\|_2^2\leqslant 1$ (note that we naturally have that $\tau(\mu,L,\gamma,n)\leqslant (\tau(\mu,L,\gamma,1))^n$). This is what we do in the following lines using \textsc{PEPit}.

\paragraph{Imports.} Before going into the example, we have to include the right \textsc{python} imports. For this example, it is necessary to perform two imports.
\begin{lstlisting}
from PEPit import PEP
from PEPit.functions import\
    SmoothStronglyConvexFunction
\end{lstlisting}

\paragraph{Initialization of \textsc{PEPit}.} First, we set the stage by initializing a PEP object. This object allows manipulating the forthcoming ingredients of the PEP, such as functions and iterates.
\begin{lstlisting}[firstnumber=last]
problem = PEP()
\end{lstlisting}
For the sake of the example, let us pick some simple values for the problem class and algorithmic parameters, for which we perform the worst-case analysis.
\begin{lstlisting}[firstnumber=last]
L = 1.          # Smoothness parameter
mu = .1         # Strong convexity parameter
gamma = 1. / L  # Step size
n = 1           # Number of iterations
\end{lstlisting}
\paragraph{Specifying the problem class.} Second, we specify our working assumptions on the function to be optimized and instantiate a corresponding object. Here, the minimization problem at hand was of the form~\eqref{eq:min_problem} with a smooth strongly convex function.
\begin{lstlisting}[firstnumber=last]
# Declare an L-smooth mu-strongly convex function
# named "func"
func = problem.declare_function(
    SmoothStronglyConvexFunction,
    mu=mu,  # Strong convexity param.
    L=L)    # Smoothness param.
\end{lstlisting}

\paragraph{Algorithm initialization.} Third, we can instantiate the starting points for the two gradient methods that we will run, and specify an \emph{initial condition} on those points. To this end, two starting points $x_0$ and $y_0$ are introduced, one for each trajectory, and a bound on the initial distance between those points is specified as $\|x_0-y_0\|^2\leqslant 1$.
\begin{lstlisting}[firstnumber=last]
# Declare two starting points
x_0 = problem.set_initial_point()
y_0 = problem.set_initial_point()

# Initial condition ||x_0 - y_0||^2 <= 1
problem.set_initial_condition((x_0 - y_0) ** 2 <= 1)
\end{lstlisting}

\paragraph{Algorithm implementation.} In this fourth step, we specify the algorithm in a natural format. In this example, we simply use the iterates (which are \textsc{PEPit} objects) as if we had to implement gradient descent in practice using a simple loop.
\begin{lstlisting}[firstnumber=last]
# Initialize the algorithm
x = x_0
y = y_0
# Run n steps of the GD method for the two sequences
for _ in range(n):
    # Replace x and y with their next iterate
    x = x - gamma * func.gradient(x)  # call to f'(x)
    y = y - gamma * func.gradient(y)  # call to f'(y)
\end{lstlisting}

\paragraph{Setting up a performance measure.} It is crucial for the worst-case analysis to specify the \emph{metric} for which we wish to compute a worst-case performance. In this example, we wish to compute the worst-case value of $\|x_n-y_n\|^2$, which we specify as follows.
\begin{lstlisting}[firstnumber=last]
# Set the performance metric to the distance
# ||x_n - y_n||^2 
problem.set_performance_metric((x-y)**2)
\end{lstlisting}

\paragraph{Solving the PEP.} The last natural stage in the process is to solve the corresponding PEP. This is done via the following line, which will ask \textsc{PEPit} to perform the modeling steps and to call an appropriate SDP solver (which should be installed beforehand) to perform the worst-case analysis.
\begin{lstlisting}[firstnumber=last]
# Solve the PEP
pepit_tau = problem.solve()
\end{lstlisting}

\paragraph{Output.} Running these pieces of code (see \href{https://pepit.readthedocs.io/en/latest/examples/k.html#contraction-rate-of-gradient-descent}{PEPit/examples/} for the complete example) for some specific values of the parameters $n=1$, $L=1$, $\mu=.1$ and $\gamma=1$, one obtains the following output.
{\small
\begin{lstlisting}[style=style1,language=bash]
(PEPit) Setting up the problem: size of the Gram matrix: 4x4
(PEPit) Setting up the problem: performance measure is the minimum of 1 element
(PEPit) Setting up the problem: Adding initial conditions
                                and general constraints ...
(PEPit) Setting up the problem: initial conditions and general constraints
                                (1 constraint(s) added)
(PEPit) Setting up the problem: interpolation conditions for 1 function(s)
            Function 1 : Adding 2 scalar constraint(s) ...
            Function 1 : 2 scalar constraint(s) added
(PEPit) Setting up the problem: additional constraints for 0 function(s)
(PEPit) Compiling SDP
(PEPit) Calling SDP solver
(PEPit) Solver status: optimal (wrapper:cvxpy, solver: MOSEK);
                       optimal value: 0.8100000029203449
(PEPit) Primal feasibility check:
        The solver found a Gram matrix that is positive semi-definite
        up to an error of 1.896548260018477e-09
        All the primal scalar constraints are verified
        up to an error of 3.042855638898251e-09
(PEPit) Dual feasibility check:
        The solver found a residual matrix that is positive semi-definite
        All the dual scalar values associated with inequality constraints
        are nonnegative
(PEPit) The worst-case guarantee proof is perfectly reconstituted
        up to an error of 4.0078754315331366e-08
(PEPit) Final upper bound (dual): 0.8100000036427537 
        and lower bound (primal example): 0.8100000029203449 
(PEPit) Duality gap: absolute: 7.224087994472939e-10
                     and relative: 8.918627121515396e-10
\end{lstlisting}}
Note that the size of the SDP is larger than that of Section~\ref{s:example} ($4\times 4$ instead of $3\times 3$ in~\eqref{eq:def_G}) because the modeling step is done in a slightly more generic way, which might not be exploiting all specificities of the problem at hand (see formulation {in~\eqref{eq:def_G} where we use the variable $x_0 - y_0$ instead of both $x_0$ and $y_0$ as they only appear together in the original problem formulation}). For more complete examples of worst-case analyses using \textsc{PEPit}, see Section~\ref{sec:numerical-examples}.

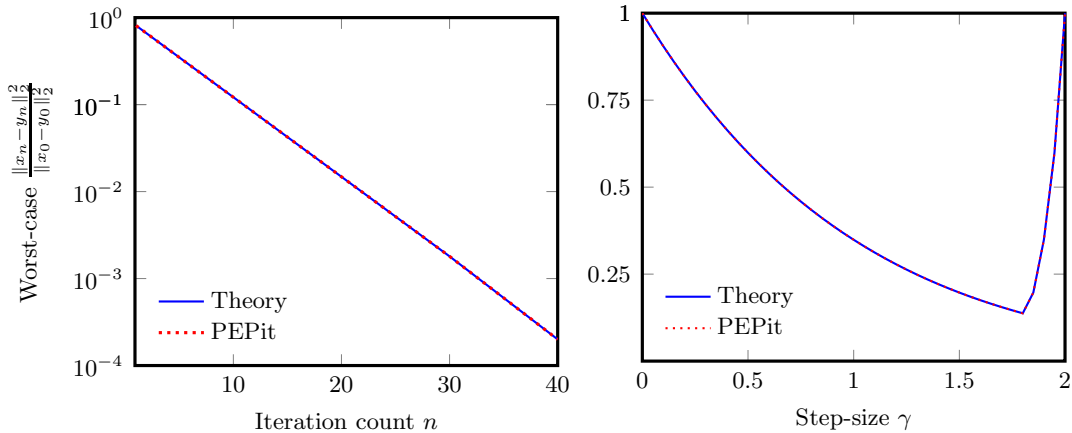
\begin{figure}[!h]
\begin{subfigure}[t]{0.45\textwidth}
  \centering
  \begin{tikzpicture}
\begin{semilogyaxis}[legend pos=south west,legend style={draw=none},legend cell align={left},plotOptions1,width=1.05\linewidth,xtick={0,10,20,30,40},xmax=40,ylabel={Worst-case $\frac{\|x_n-y_n\|^2_2}{\|x_0-y_0\|^2_2}$},xlabel={Iteration count $n$}]
\addplot[blue, thick] table [y=theorygdc, x=condition]{data/gdc_n.txt}; 
\addplot[dashed, red, very thick] table [y=gdc, x=condition]{data/gdc_n.txt};
\end{semilogyaxis}
\end{tikzpicture}
  \caption{Worst-case guarantee and theoretical tight bound as a function of the iteration count $n$ for $\gamma=\frac{1}{L}$.}
  \label{fig:GDC_n}
\end{subfigure}
\hspace{.5cm}
\begin{subfigure}[t]{0.45\textwidth}
  \centering
  \begin{tikzpicture}
\begin{axis}[legend pos=south west,legend style={draw=none},legend cell align={left},plotOptions5,width=1.05\linewidth,ylabel={},xlabel={Step-size $\gamma$}]
\addplot[blue, thick] table [y=theorygdc, x=condition]{data/gdc_step.txt};
\addplot[red, dashed, very thick] table [y=gdc, x=condition]{data/gdc_step.txt};
\end{axis}
\end{tikzpicture}
  \caption{Worst-case guarantee and theoretical tight bound as a function of the step-size $\gamma$, at iteration $n=5$.}
  \label{fig:GDC_step}
  \end{subfigure}
 \caption{Comparison: worst-case guarantee from \textsc{PEPit}~(\textbf{plain blue}) and theoretical tight worst-case bound~\eqref{eq:tau_GD_ex} for gradient descent in terms of $\frac{\|x_n-y_n\|^2_2}{\|x_0-y_0\|^2_2}$~(\textbf{dashed red}). Problem parameters fixed to $\mu=0.1$ and~$L=1$.}
\label{fig:examplesgd}
\end{figure}

It is also possible to run the code for different values of the parameters, as exemplified in Figure~\ref{fig:examplesgd}. This simple example allows us to observe that numerical values obtained from \textsc{PEPit} match the worst-case guarantee~\eqref{fig:GDC_n}, and to optimize the step-size numerically in Figure~\ref{fig:GDC_step}.

\FloatBarrier

\section{\textsc{PEPit} code structure and semidefinite formulation}\label{s:SDP_details}

This section provides the general semidefinite program (SDP) that is formulated by the package, as well as its relationship with the code. As already underlined, \textsc{PEPit} aims at providing simple ways to model performance estimation problems~(PEPs) by abstracting the coding of those SDPs. Then, \textsc{PEPit} passes the SDP either (i) to CVXPY~\cite{diamond2016cvxpy}, thereby benefiting from all SDP solvers that are interfaced with it (such as~\cite{o2016conic,mosek}), or (ii)~directly to MOSEK~\cite{mosek}.

Before going into the implementation details of PEPit, note that executing the codes from Section 3.2 to 3.5 requires the following imports.
\begin{lstlisting}
# For Section 3.2 to Section 3.5
import PEPit
import PEPit.functions

# For Section 3.3 and Section 3.5 we also need:
from PEPit.functions import\
    ConvexFunction,\
    ConvexIndicatorFunction,\
    SmoothConvexFunction

# For the examples of Section 3.3:
from PEPit.primitive_steps import\
    proximal_step,\
    linear_optimization_step,\
    inexact_gradient_step
\end{lstlisting}

\subsection{Semidefinite formulation}\label{s:sdp_form}

The package formulates SDPs of the form
\begin{equation}\label{eq:gen_SDP}
\begin{aligned}
    \max_{\substack{\tau\in\mathbb{R},G\in\mathbb{S}^{n_p},H\in\mathbb{R}^{n_h}}} &\quad \tau \\
        \text{s.t.} \,\,& G\succcurlyeq 0\\
        & \tau -\tr{A_i G}-a_i^T H-\alpha_i\leqslant 0 \quad \forall i\in I_1\\
        & \left[\tr{B_{j,k,i}G}+b_{j,k,i}^T H+\beta_{j,k,i}\right]_{1\leq j,k \leq n_i}\succcurlyeq 0 \quad \forall i\in I_2,\\
\end{aligned}
\end{equation}
with $n_h,\, n_p\in\mathbb{N}$, some index sets $I_1$ and  $I_2$ in $\mathbb{N}$, and sets of problem parameters $\{(A_i,a_i,\alpha_i)\}_{i\in I_1}$, and $\{(B_{j,k,i},b_{j,k,i},\beta_{j,k,i})\}_{1\leq j,k\leq n_i, i\in I_2}$ ($n_i\in\mathbb{N}_{>0}$ for all $i\in I_2$) of appropriate dimensions, which are constructed from the algorithm and the class of problems at hand (and hence all depend on the parameters of the algorithm and of those of the class of problems) for computing appropriate worst-case scenarios.

Formulating such an SDP is usually cumbersome and relatively error-prone, and the role of \textsc{PEPit} is to generate all those parameters in a user-friendly way. {That is, \textsc{PEPit} parses the natural description of the problem the user is familiar with.} The problem is not described in terms of the SDP variables $(G, H)$ but rather in terms of another couple $(P, H)$ with the following structure:
\begin{equation}
P \eqdef [p_1,\, p_2,\, \ldots,\, p_{n_p}], \quad\quad\quad H \eqdef [h_1,\,h_2\,\,\ldots,\,h_{n_h}],
\end{equation}
for some {$\{p_i\}_{1\leq i\leq n_p}\subset\mathbb{R}^d$} for some $d\in\mathbb{N}$, and $\{h_i\}_{1\leq i\leq n_h}\subset\mathbb{R}$. Hence $P\in \mathbb{R}^{d\times n_p}$ and $H\in\mathbb{R}^{n_h}$. The relationship with~\eqref{eq:gen_SDP} is that $G$ can be constructed without loss of generality as $G\eqdef P^T P\succeq 0$.

\subsection{Base \textsc{PEPit} objects}
Using previous notations, \textsc{PEPit} allows the user to formulate the problem in terms of $(P, H)$ instead of $(G, H)$, which turns out much more natural for describing many algorithms and problem classes. {The base working procedure is as follows}:
\begin{itemize}
    \item {E}ach $p_i$ corresponds to a base \textsc{PEPit.Point} object (referred to as a \emph{leaf} element). Such objects can be added and subtracted together, and scaled by real values for forming new objects \textsc{PEPit.Point} (which are then combinations of \emph{leaf} elements). \textsc{PEPit.Point} objects can be understood as elements of $\mathbb{R}^d$, the space of iterates/gradients.
	\item Each $h_i$ corresponds to a base \textsc{PEPit.Expression} object (referred to as a \emph{leaf} element). Such objects can be added and subtracted together, but also scaled by real values for forming new objects \textsc{PEPit.Expression} (which are combinations of \emph{leaf} elements). \textsc{PEPit.Expression} objects can be understood as elements of $\mathbb{R}$, such as, for instance, (scalar) function values.
	\item It is possible to compute dot products of \textsc{PEPit.Point} objects (e.g., $p_i^T p_j$ for some $1\leq i,j\leq n_p$), resulting in \textsc{PEPit.Expression} objects.
	\item Comparing two \textsc{PEPit.Expression} objects with an operator in $\{=,\leq,\geq\}$ leads to a \textsc{PEPit.Constraint} object. Similarly, \textsc{PEPit.Expression} can be gathered in arrays to form \textsc{PEPit.Psd\_matrix} objects for formulating semidefinite constraints.
\end{itemize}

\begin{example} Following up on the notations from Section~\ref{s:example} for describing one iteration of gradient descent, a possibility is to think of $p_1$ as corresponding to some $x_0$, of $p_2$ as a corresponding gradient $g_{x_0}$, and of $h_1$ as the corresponding function value $f_{x_0}$. For a unit step-size, one can form $x_1=p_1-p_2$ as follows.
\begin{lstlisting}
x_0 = PEPit.Point()  # a leaf PEPit.Point (p_1)
g_x_0 = PEPit.Point()  # a leaf PEPit.Point (p_2)
# a leaf PEPit.Expression (h_1):
f_x_0 = PEPit.Expression()

x_1 = x_0 - g_x_0  # x_1 is a PEPit.Point
g_x_1 = PEPit.Point()  # a leaf PEPit.Point (p_3)
# a leaf PEPit.Expression (h_2):
f_x_1 = PEPit.Expression()
\end{lstlisting}
It is important to note that each call to \textsc{PEPit.Point()} increments $n_p$ (dimension $G$ in~\eqref{eq:gen_SDP}). Similarly, \textsc{PEPit.Expression()} increments $n_h$. It is cheaper to solve a problem with as few leaf points and leaf expressions as possible.

{To formulate} the objective, initial conditions, and interpolation constraints, we use \textsc{PEPit.Expression} objects, which we compare {together to form}
\newline \textsc{PEPit.Constraint} objects.
For instance, interpolation conditions for convex functions can be formulated as objects of the \textsc{PEPit.Constraint} class:
\begin{lstlisting}[firstnumber=last]
# this is a PEPit.Expression object:
expr = f_x_0 + g_x_0 * (x_1 - x_0)
# those are two PEPit.Constraint objects:
cons_1 = (f_x_0 + g_x_0 * (x_1 - x_0) <= f_x_1)
cons_2 = (f_x_1 + g_x_1 * (x_0 - x_1) <= f_x_0)
\end{lstlisting}
\end{example}

{Specifying constraints as done above, in a natural mathematical way, is convenient for the user.}
Once all required points, expressions, and constraints are associated with a PEP (as in the example of Section~\ref{s:code} or in the following lines), \textsc{PEPit} takes care of formulating the appropriate index sets $I_1$, and $I_2$ as well as all problem parameters for~\eqref{eq:gen_SDP}: 
$\{(A_i,a_i,\alpha_i)\}_{i\in I_1}$, and $\{(B_{j,k,i},b_{j,k,i},\beta_{j,k,i})\}_{1\leq j,k\leq n_i, i\in I_2}$ for feeding~\eqref{eq:gen_SDP} to SDP solvers~\cite{diamond2016cvxpy,o2016conic,mosek}. 
{However, on the user side, specifying a large number of such constraints remains relatively cumbersome. Therefore, \textsc{PEPit} relies on a few additional structures and aliases that allow abstractly generating constraints blocks.}

\subsection{Main \textsc{PEPit} simplifying abstractions and aliases}

Two key abstractions appearing in \textsc{PEPit} are the \emph{functions} (which can also be manipulated algebraically) and \emph{oracles} (or \emph{primitive steps}, which are simple routines). {In both cases, those structures were designed to allow users to easily add new types of functions and oracles.}

\paragraph{\textsc{PEPit} functions.} Among the most important building blocks simplifying the formulation of the constraints in~\eqref{eq:gen_SDP}, \textsc{PEPit.Function} objects are particularly important. They allow for generating large numbers of constraints by remaining close to clean mathematical statements. Their most important characteristics are as follows:
\begin{itemize}
	\item \textsc{PEPit.Function()} creates a new \emph{leaf} function,
	\item {E}ach leaf function contains a list of triplet $\{(x_i,g_i,f_i)\}_i$ which corresponds to the sampled version of the function in the form (points, gradients, function values), as presented in Section~\ref{s:example}.
	\item Each \textsc{PEPit.Function} object \textsc{f} is featured with a \textsc{f.add\_point(triplet)} method, which adds a triplet to the list of samples associated to~\textsc{f}. By relying on appropriate abstractions (exemplified later), the users are expected to almost never use this method explicitly.
	\item Each \textsc{PEPit.Function} object \textsc{f} contains a \textsc{f.add\_class\_constraints()} method that generate the list of interpolation \textsc{PEPit.Constraint} associated to the sampled version of \textsc{f}. This method is never explicitly used by the user but allows \textsc{PEPit} to translate the mathematical statements (description of the function) to actual constraints.
    \item Each \textsc{PEPit.Function} can also be scaled by scalar values and added or subtracted to other \textsc{PEPit.Function} for creating new \textsc{PEPit.Function} objects.
\end{itemize} 
{Functions also feature a number of aliases} that allow to easily create specific \textsc{PEPit.Point} and \textsc{PEPit.Expression} objects associated with the function (such as accessing the gradient/value at a specific point, or accessing an optimal point of a function).

\begin{example}[Base operations with functions] This example shows that we can manipulate functions to create new functions.
\begin{lstlisting}
f = PEPit.Function()  # a leaf PEPit.Function
h = PEPit.Function()  # a leaf PEPit.Function
F = 2 * f + h  # a PEPit.Function
\end{lstlisting}
We can also call base function operations on functions that are constructed by combining \emph{leaf} functions.
\begin{lstlisting}[firstnumber=last]
x = PEPit.Point()  # this is p_1
g_x = PEPit.Point()  # this is p_2
f_x = PEPit.Expression()  # this is h_1

F.add_point((x, g_x, f_x))
# internally creates p_3 (a new point)
# and h_2 (a new expression),
# and adds one point to the list of f,
# and one to the list of h so that
# the weighted sums of gradients and 
# function values at x are correct.
\end{lstlisting}
However, {to simplify} the usage of \textsc{PEPit}, one should essentially avoid directly using the method \textsc{add\_point}, and rather rely on more readable statements.
In that regard, the following piece of code is equivalent to the previous one but relies on more readable operations.
\begin{lstlisting}
f = PEPit.Function()  # a leaf PEPit.Function
h = PEPit.Function()  # a leaf PEPit.Function
F = 2 * f + h  # a PEPit.Function

x = PEPit.Point()  # this is p_1 (arbitrary point)

g_x = F.gradient(x)  # g_x is the gradient of F at x
f_x = F(x)  # f_x is F(x)
\end{lstlisting}
\end{example}
{The list of points the function (or its gradient) is evaluated on is stored in the underlying object and finally used to generate the corresponding list of interpolation constraints.} For creating a class that features a new type of interpolation conditions, it suffices to create a new class that inherits \textsc{PEPit.Function} and implements the operation \textsc{add\_class\_constraints} with the appropriate interpolation constraints (see \url{https://github.com/PerformanceEstimation/PEPit/tree/master/PEPit/functions} for examples).

Let us now mention an important class of convenient abstraction that belongs to the heart of \textsc{PEPit}'s philosophy.

\paragraph{\textsc{PEPit} primitive steps (or oracles).} \textsc{PEPit} \emph{primitive steps} (or \emph{oracles}) are essentially simple aliases defining notational shortcuts for some algorithmic operations, rendering them closer to their mathematical abstractions.
They generally consist of the appropriate creations (and constraining) of points and expressions, for instance by adding appropriate triplets to sampled versions of the functions under consideration.
Let us provide a few examples.

\begin{example}[Proximal operators.] In this example, we provide two ways to use proximal operators within \textsc{PEPit}.
The first one uses the base \textsc{PEPit} methods for doing this, and the second one relies on the more readable \textsc{proximal\_step} that is provided as a primitive step of \textsc{PEPit}.
Those two ways are computationally equivalent, though different in terms of readability.
Let us recall that the proximal operation (with unit step-size) associated to~$f$ is given by
\[ x_1=\mathrm{prox}_{f}(x_0)\eqdef \mathrm{argmin}_x\left\{f(x)+\tfrac{1}{2}\|x-x_0\|^2\right\}.\]
For any (closed, proper) convex function $f$, this operation is well-defined  and amounts to an \emph{implicit subgradient} operation:
\[ x_1=x_0-g_1,\]
with $g_1\in\partial f(x_1)$.
\begin{lstlisting}
# f is a (closed, proper) convex function
f = PEPit.functions.ConvexFunction()
x0 = PEPit.Point()  # is a point

# how to construct x_1 = prox_f (x_0)?
# (i) initiate some g1 and f1
g1 = PEPit.Point()
f1 = PEPit.Expression()

# (ii) form x1 using x0 and g1
x1 = x0 - g1

# (iii) constrain g1 to be a subgradient of f at x1
f.add_point((x1, g1, f1))
\end{lstlisting}
Equivalently, using \textsc{proximal\_step}, we have:
\begin{lstlisting}
# f is a (closed, proper) convex function
f = PEPit.functions.ConvexFunction()
x0 = PEPit.Point()  # is a point

x1, g1, f1 = proximal_step(x0, f, gamma=1)
# note: gamma is a step-size, set to 1 in the example.
\end{lstlisting}
\end{example}

\begin{example}[Linear optimization oracles.] This operation is at the core of the Frank-Wolfe (a.k.a.\ conditional gradient) method, and consists in, given a (closed, convex) domain $\mathcal{K}$ (whose indicator function is denoted by $i_{\mathcal{K}}$) and a search direction $d_0$, in computing a solution to
\[x_1 \in \mathrm{argmin}_x \, d_0^T x + i_{\mathcal{K}}(x),\]
which is mathematically equivalent to writing $-d_0\in\partial i_{\mathcal{K}}(x_1)$.
In \textsc{PEPit}, this can also easily be coded as follows.
\begin{lstlisting}
# ind is a (closed) convex indicator function
ind = PEPit.functions.ConvexIndicatorFunction()
d0 = PEPit.Point()

# how to construct a solution to min_x d0*x + ind(x) ?
# (i) Initiate a new point and an expression
x1 = PEPit.Point()
f1 = PEPit.Expression()

# (ii) Constrain -d0 to be a subgradient
# of the indicator at x1
ind.add_point((x1, -d0, f1))
\end{lstlisting}
Using abstraction again, this code is equivalent to
\begin{lstlisting}
# ind is a (closed) convex indicator function
ind = PEPit.functions.ConvexIndicatorFunction()
d0 = PEPit.Point()

# how to construct a solution to min_x d0*x + ind(x) ?
x1, _, f1 = linear_optimization_step(d0, ind)
\end{lstlisting}
\end{example}

\begin{example}[Using approximate gradients.] Another standard operation in first-order optimization consists in using \emph{approximate} gradient values.
For instance, for computing some $x_1$ using a gradient iteration (with unit step-size) with an approximate gradient $\tilde{d}_0\approx \nabla f(x_0)$ for some appropriate function $f$:
\[ x_1 = x_0 - \tilde{d}_0,\]
with $\tilde{d}_0$ being an $\epsilon$ approximation to $\nabla f(x_0)$ in the following sense: $\|\tilde{d}_0-\nabla f(x_0)\|\leq \epsilon \|\nabla f(x_0)\|$, say with $\epsilon=0.1$.
\begin{lstlisting}
# f is a 1-smooth convex function
f = PEPit.functions.SmoothConvexFunction(L=1)
x0 = PEPit.Point()
epsilon = .1

# how to construct a x1?
# (i) initiate a d0 and the gradient of x0
g0 = f.gradient(x0)
d0 = PEPit.Point()
f.add_constraint((d0 - g0) ** 2 <= epsilon * g0 ** 2)

x1 = x0 - d0
\end{lstlisting}
This can equivalently be done using \textsc{inexact\_gradient\_step}, as follows.
\begin{lstlisting}
# f is a 1-smooth convex function
f = PEPit.functions.SmoothConvexFunction(L=1)
x0 = PEPit.Point()

# how to construct a x1?
from PEPit.primitive_steps\
    import inexact_gradient_step

x1, d0, _ = inexact_gradient_step(x0, f, gamma=1,
                                  epsilon=.1,
                                  notion='relative')
\end{lstlisting}

\end{example}

To conclude this section, \textsc{PEPit}'s philosophy is to contain many \emph{abstract} routines that can be associated with simple mathematical statements. Taken separately, those routines are relatively simple and can easily be created or modified by the users.

\subsection{The objective function of the PEP: performance metrics}

A key point that we did not mention so far concerns the objective function of~\eqref{eq:gen_SDP}.
It is handled by what is referred to as \emph{performance metrics} in \textsc{PEPit}.
In the SDP formulation~\eqref{eq:gen_SDP}, they correspond to the index set $I_1$ and the set of parameters $\{(A_i,a_i,\alpha_i)\}_{i\in I_1}$, and they are handled by calling the \textsc{set\_performance\_metric} method (which takes a \textsc{PEPit.Expression} as sole argument) associate to a PEP object.
By introducing the variable $\tau$ in~\eqref{eq:gen_SDP}, the objective of the PEP corresponds to the minimum value of all specified performance metrics.

\begin{example} This example shows how to specify an objective function (a.k.a. \emph{performance metric}) within \textsc{PEPit}.
\begin{lstlisting}
problem = PEPit.PEP()

p1 = PEPit.Point()
p2 = PEPit.Point()
h1 = PEPit.Expression()

# set the PEP objective as the minimum value among
# ||p1||^2, ||p2||^2, and h1.
problem.set_performance_metric(p1**2)
problem.set_performance_metric(p2**2)
problem.set_performance_metric(h1)
\end{lstlisting}
\end{example}
\subsection{Formulating and solving the PEP}

As a final stage for formulating~\eqref{eq:gen_SDP}, we need to gather all functions and constraints together and reformulate in terms of~\eqref{eq:gen_SDP}.
\textsc{PEPit} performs it through the use of the PEP object (as provided in Section~\ref{s:code}) as follows:
\begin{lstlisting}
problem = PEPit.PEP()
\end{lstlisting}
which is used for centralizing all the information about the PEP at hand.
In particular, it is a good practice to avoid using \textsc{PEPit.Point()} directly, and to rather use \textsc{problem.set\_initial\_point()}.
Similarly, new functions should be declared through the PEP object using the \textsc{declare\_function} method:
\begin{lstlisting}
# declares a convex (closed, proper) function h
h = problem.declare_function(ConvexFunction) 

# declares a convex (closed) indicator ind
ind = problem.declare_function(ConvexIndicatorFunction) 

# declares a 1-smooth convex function f
f = problem.declare_function(SmoothConvexFunction,
                             L=1)
\end{lstlisting}

{To model~\eqref{eq:gen_SDP}}, we gather all \textsc{PEPit.Constraint} objects that are associated with \textsc{problem} through the different abstractions used in the code. For instance, the \textsc{PEPit.PEP} object \textsc{problem} will call the \textsc{add\_class\_constraints} method of all functions involved in the PEP. Once this is done, \textsc{PEPit} generates all appropriate matrices and index sets {to model} the problem as an SDP in the form~\eqref{eq:gen_SDP}, and passes it to either CVXPY~\cite{diamond2016cvxpy} or directly to MOSEK~\cite{mosek}, before performing a few post-processing steps.

\subsection{Post-processing}

Once the PEP is solved numerically, we need to manipulate its output, either for constructing proofs (on the dual PEP side), or counter-examples (on the primal PEP side).

\paragraph{Dual reconstructions.}
Dual values associated with the different constraints play an important role, as they allow us to construct mathematical proofs.

In the following, we assume that \emph{func} has been defined as a Function used in the PEP, that \emph{constraints\_list} has been defined as a list of Constraint objects involved in the PEP, and that the PEP has been solved through the \textsc{PEP.solve()} method.

One of the important features of PEPit is to {automatically} deal with interpolation constraints for the different functions involved in the PEP. This way, the attribute \emph{list\_of\_class\_constraints} enables accessing the list of Constraint objects encoding the interpolation constraints of \emph{func}.

{Each dual value} can be accessed through the \textsc{eval\_dual()} method of the associated Constraint object:
\begin{lstlisting}
for constraint in func.list_of_class_constraints:
    dual_value = constraint.eval_dual()
    print(dual_value)
\end{lstlisting}
Moreover, to ease the reconstruction of the proof, and avoid mistakes by associating a constraint to the wrong dual value, a user can name a constraint through the \textsc{set\_name()} method and access it later on through the \textsc{get\_name()} method.
Note it is the responsibility of the user to set a name to the constraints to be able to recover it later.
\begin{lstlisting}
for constraint in constraints_list:
    constraint_name = constraint.get_name()
    dual_value = constraint.eval_dual()
    print(constraint_name, dual_value)
\end{lstlisting}

Since PEPit deals with interpolation constraints without any intervention from the user side, a short description of each of those constraints is set by default, based on {Points' names} and \emph{func}'s name which can also be set through a \textsc{set\_name} method.
Finally, a user can also obtain all the dual values associated with the interpolation constraints of a function at once, using the \textsc{get\_class\_constraints\_duals} method as
\begin{lstlisting}
# assuming func has been defined
# as a Function object involved in the PEP
# and that the PEP has been solved.

tables = func.get_class_constraints_duals()
\end{lstlisting}
This method returns a dictionary whose values are pandas.DataFrames, offering great readability.
Completing the example of Section~\ref{s:code} (this time for 2 steps) with this feature, we obtain the following code:
\begin{lstlisting}
import numpy as np

from PEPit import PEP
from PEPit.functions import\
    SmoothStronglyConvexFunction

# We set the parameter of the problem
# Here we study the contraction of 1 step
# of the GD method with step 1/L on the class of
# L=1 smooth and mu=.1 strongly convex functions.

L = 1.          # Smoothness parameter
mu = .1         # Strong convexity parameter
gamma = 1. / L  # Step size
n = 2           # Number of iterations

# Instantiate the PEP object
problem = PEP()

# Declare an L-smooth mu-strongly convex function
# named "func"
func = problem.declare_function(
    SmoothStronglyConvexFunction,
    mu=mu,     # Strong convexity param.
    L=L,       # Smoothness param.
    name="f")  # Name

# Declare two starting points
x_0 = problem.set_initial_point(name="x_0")
y_0 = problem.set_initial_point(name="y_0")

# Initial condition ||x_0 - y_0||^2 <= 1
problem.set_initial_condition((x_0 - y_0) ** 2 <= 1)

# Initialize the algorithm
x = x_0
y = y_0
# Run n steps of the GD method for the two sequences
for i in range(n):

    # Replace x and y with their next iterates
    x = x - gamma * func.gradient(x)  # call to f'(x)
    x.set_name("x_{}".format(i+1))

    y = y - gamma * func.gradient(y)  # call to f'(y)
    y.set_name("y_{}".format(i+1))

# Set the performance metric to the distance
# ||x_n - y_n||^2
problem.set_performance_metric((x-y)**2)

# Solve the PEP
pepit_tau = problem.solve()

# By linearly combining the interpolation constraints
# with the right coefficients, we can prove
# ||x_n - y_n||^2 <= pepit_tau ||x_0 - y_0||^2
# The coefficient we need are the dual values
# of the interpolation constraints of func.
tables = func.get_class_constraints_duals()

# A user can display the dictionary as is,
# or can access one specific table by their name.
# Those names are intuitive, yet to be known,
# for example by displaying the keys of tables.
# Here we use the only key of this dictionary.
table = tables["smoothness_strong_convexity"]

print("\nDual values associated with"
      " interpolation constraints:")
print(table.astype(dtype=np.float16))
\end{lstlisting}

Running this code outputs the following message:
{\small
\begin{lstlisting}[style=style1,language=bash]
(PEPit) Setting up the problem: size of the Gram matrix: 6x6
(PEPit) Setting up the problem: performance measure is the minimum of 1 element
(PEPit) Setting up the problem: Adding initial conditions
                                and general constraints ...
(PEPit) Setting up the problem: initial conditions and general constraints
                                (1 constraint(s) added)
(PEPit) Setting up the problem: interpolation conditions for 1 function(s)
			Function 1 : Adding 12 scalar constraint(s) ...
			Function 1 : 12 scalar constraint(s) added
(PEPit) Setting up the problem: additional constraints for 0 function(s)
(PEPit) Compiling SDP
(PEPit) Calling SDP solver
(PEPit) Solver status: optimal (wrapper:cvxpy, solver: MOSEK);
                       optimal value: 0.6561000087150534
(PEPit) Primal feasibility check:
	The solver found a Gram matrix that is positive semi-definite
        up to an error of 2.584803430062655e-09
	All the primal scalar constraints are verified
        up to an error of 4.590572921792102e-09
(PEPit) Dual feasibility check:
	The solver found a residual matrix that is positive semi-definite
	All the dual scalar values associated with inequality constraints
        are nonnegative up to an error of 3.8932973849815053e-10
(PEPit) The worst-case guarantee proof is perfectly reconstituted
        up to an error of 1.349705540942825e-07
(PEPit) Final upper bound (dual): 0.6561000067100656
        and lower bound (primal example): 0.6561000087150534 
(PEPit) Duality gap: absolute: -2.004987731396568e-09
                     and relative: -3.055917855150253e-09

Dual values associated with interpolation constraints:
IC_f       x_0       y_0       x_1       y_1
x_0   0.000000  1.458008  0.000000  0.000000
y_0   1.458008  0.000000  0.000000  0.000000
x_1  -0.000000  0.000000  0.000000  1.799805
y_1   0.000000 -0.000000  1.799805  0.000000
\end{lstlisting}
}

\paragraph{Primal reconstructions.} In order to construct an example of a problem on which the algorithm behaves {as poorly as possible with respect to the given performance metrics}, all \textsc{PEPit.Point}, all \textsc{PEPit.Expression} and all \newline \textsc{PEPit.Constraint} objects can be conveniently evaluated through the \textsc{eval()} method.
Moreover, their \textsc{name} attribute can help to sort them. \textsc{PEPit} also offers the possibility to search for simpler, potentially low-dimensional problem instances via the trace norm~\cite{recht2010guaranteed} or the logdet~\cite{fazel2003log} heuristics (aiming at finding low-rank feasible matrices $G$ for the problem~\eqref{eq:gen_SDP} while keeping the same objective value).
Those post-processing steps can be accessed via the option of the \textsc{solve} method (see \url{https://pepit.readthedocs.io/en/latest/api/main_modules.html#pep}). Examples of such usages can be found in the documentation at \url{https://pepit.readthedocs.io/en/latest/examples/j.html}.

\FloatBarrier

\section{\textsc{PEPit}: general overview and content}\label{s:content}
In this section, we go back to the mathematical content of the toolbox and describe the various choices of (i)~elementary oracles used in algorithms, (ii)~problem or function classes, (iii)~performance measures, and (iv)~initial conditions, that are naturally handled by \textsc{PEPit}.
\textsc{PEPit} also allows studying methods for monotone inclusions and fixed point problems, but we do not cover them in this summary.
In the optimization setting, the minimization problem under consideration has the form
\begin{equation}\label{eq:min_composite}
    F_\star\eqdef \min_{x\in\mathbb{R}^d} \left\{F(x)\equiv \sum_{i=1}^K f_i(x)\right\},
\end{equation}
for some $K\in\mathbb{N}$ and where each $f_i$ is assumed to belong to some class of functions denoted by $\mathcal{F}_i$, encoding our working assumptions on $f_i$.
We further assume the algorithms to gather information about the functions $\{f_i\}_i$ only via black-box oracles such as gradient or proximal evaluations.

\paragraph{Black-box oracles.} The base black-box optimization oracles/operations available in \textsc{PEPit} are the following:
\begin{itemize}
    \item (sub)gradient steps,
    \item proximal and projection steps,
    \item linear optimization (or conditional) steps.
\end{itemize}
\textsc{PEPit} also allows for their slightly more general approximate/inexact and Bregman (or mirror) versions.
Those oracles might be combined with a few other operations, such as exact line-search procedures, as detailed in~Table~\ref{tab:oracles}.
\begin{table*}[ht]
\begin{center}
\caption{Main base primitive steps (oracles) included in \textsc{PEPit}. Appropriate references are provided in the corresponding \href{https://pepit.readthedocs.io/en/latest/api/steps.html}{documentation}.
Some oracles are overlapping and are present for promoting a better readability of the code and for numerical efficiency. Variations around the present oracles can be added at will. For each oracle, $x_+$ denotes the output of the oracles; the other elements are either input of the oracles or intermediary optimization variables.}\label{tab:oracles}
{\renewcommand{\arraystretch}{1.8}
\begin{tabularx}{\textwidth}{p{4cm} p{5.5cm}  p{4cm}}
\toprule
Oracle name & Description & Tightness\\
\midrule
Gradient step & {$x_+ = x - \gamma g$ with $g=\nabla f(x)$} & \ding{52} \\
Subgradient step & $x_+ = x - \gamma g$ with $g\in\partial f(x)$ & \ding{52} \\
\href{https://pepit.readthedocs.io/en/latest/api/steps.html\#epsilon-subgradient-step}{Epsilon-subgradient step} & $x_+ = x - \gamma g$ with $g\in\partial_{\varepsilon} f(x)$ & \ding{52} \\
\href{https://pepit.readthedocs.io/en/latest/api/steps.html\#inexact-gradient-step}{Inexact gradient step} & $x_+ = x - \gamma g$ with $g\approx_\epsilon \nabla f(x)$ \newline for some notion ``$\approx_\epsilon$'' of approximation. & \ding{52}\\
\href{https://pepit.readthedocs.io/en/latest/api/steps.html\#exact-line-search-step}{Exact line-search step} & $x_+ = \underset{z \in x + \Span \left\{d_i, i\in \range{1}{T} \right\}}{\argmin} f(z)$ & \ding{54}\\
\href{https://pepit.readthedocs.io/en/latest/api/steps.html\#proximal-step}{Proximal step} & $x_+ = \underset{z}{\argmin} \left\{ \gamma f(z) + \frac{1}{2} \|z - x\|^2 \right\}$ & \ding{52}\\
\href{https://pepit.readthedocs.io/en/latest/api/steps.html\#inexact-proximal-step}{Inexact proximal step} & $x_+ \approx_\epsilon \underset{z}{\argmin} \left\{ \gamma f(z) + \frac{1}{2} \|z - x\|^2 \right\}$ \newline for some notion ``$\approx_\epsilon$'' of approximation. & \ding{52} \\
\href{https://pepit.readthedocs.io/en/latest/api/steps.html\#bregman-gradient-step}{Bregman gradient step} & $x_+ = \underset{z}{\argmin} \left[ \left< \nabla f(x); z - x \right> + \frac{1}{\gamma} D_h(z; x) \right]$ & \ding{52} \\
\href{https://pepit.readthedocs.io/en/latest/api/steps.html\#bregman-proximal-step}{Bregman proximal step} & $x_+ = \underset{z}{\argmin} \left[ f(z) + \frac{1}{\gamma} D_h(z; x) \right]$ & \ding{52} \\
\href{https://pepit.readthedocs.io/en/latest/api/steps.html\#linear-optimization-step}{Linear optimization step} & $x_+ = \underset{z ~|~ \text{ind}(z)=0}{\argmin} \left< \text{dir}; z \right>$ & \ding{52}\\
\bottomrule
\end{tabularx}
}
\end{center}
\end{table*}

\paragraph{Problem classes.} A few base classes of functions are readily available within the package (see~Table~\ref{tab:functions} for further details) such as:
    \begin{itemize}
        \item {C}onvex functions within different classes of assumptions possibly involving bounded gradients (Lipschitz functions), {quadratically upper bounded functions,} bounded domains, smoothness, and/or strong convexity.
        Those assumptions might be combined when compatible.
        \item Convex indicator functions, possibly with a bounded domain.
        \item Smooth nonconvex functions.
        \item Quadratic functions.
    \end{itemize}
Beyond the pure optimization setting, \textsc{PEPit} also allows using operators (see~Table~\ref{tab:operators}) within different classes of assumptions (namely: nonexpansive, Lipschitz, cocoercive, maximally monotone, and strongly monotone operators) for studying first-order methods for monotone inclusions and variational inequalities.
\begin{table*}[ht]
\begin{center}
{\renewcommand{\arraystretch}{1.8}
\caption{Some base function classes included in \textsc{PEPit}, detailed in the \href{https://pepit.readthedocs.io/en/latest/api/functions.html}{documentation}. Default functional classes within \textsc{PEPit}. Some classes are overlapping and are present only to promote a better readability of the code.\label{tab:functions}
}
\begin{tabular}{@{}lcc@{}}
\toprule
Function class name & Tightness \\
\midrule
    \href{https://pepit.readthedocs.io/en/latest/api/functions.html#convex}{Convex (closed, proper) functions} &  \ding{52} \\
	\href{https://pepit.readthedocs.io/en/latest/api/functions.html#convex-and-lipschitz-continuous}{Convex (closed, proper) Lipschitz-continuous functions} & \ding{52} \\
	\href{https://pepit.readthedocs.io/en/latest/api/functions.html#convex-indicator}{Convex (closed, proper) indicator functions} & \ding{52}\\
	\href{https://pepit.readthedocs.io/en/latest/api/functions.html#convex-support-functions}{Convex support functions} & \ding{52}\\
	\href{https://pepit.readthedocs.io/en/latest/api/functions.html#strongly-convex-and-smooth}{Smooth strongly convex functions} & \ding{52}\\
    \href{https://pepit.readthedocs.io/en/latest/api/functions.html#convex-and-smooth}{Smooth convex functions} & \ding{52}
    \\
	\href{https://pepit.readthedocs.io/en/latest/api/functions.html#smooth}{Smooth (possibly nonconvex) functions} & \ding{52}\\
    \href{https://pepit.readthedocs.io/en/latest/api/functions.html#smooth-convex-and-lipschitz-continuous}{Smooth convex Lipschitz functions} & \ding{52}\\
	\href{https://pepit.readthedocs.io/en/latest/api/functions.html#strongly-convex}{Strongly convex functions} & \ding{52}\\
    \href{https://pepit.readthedocs.io/en/latest/api/functions.html#convex-and-quadratically-upper-bounded}{Convex quadratically upper-bounded functions} & \ding{52}\\
    \href{https://pepit.readthedocs.io/en/latest/api/functions.html#restricted-secant-inequality-and-error-bound}{Restricted secant inequality and error bound} & \ding{52}\\
    \href{https://pepit.readthedocs.io/en/latest/api/functions.html#strongly-convex-and-smooth-quadratic}{Smooth strongly convex quadratic functions} & \ding{52}\\
    \href{https://pepit.readthedocs.io/en/latest/api/functions.html#convex-and-smooth-by-block}{Smooth convex function by block} & \ding{54}\\
\bottomrule
\end{tabular}}
\end{center}
\end{table*}

\begin{table*}[ht]
\begin{center}
{\renewcommand{\arraystretch}{1.8}
\caption{Base operator classes within \textsc{PEPit}, detailed in the \href{https://pepit.readthedocs.io/en/latest/api/operators.html}{documentation}.
Some classes are overlapping and are present only to promote a better readability of the code.
Note that, for some classes, the associated constraints might not be tight, meaning that the methodology might only be able to generate upper-bound on the worst-case behaviors.\label{tab:operators}}
\begin{tabular}{@{}lcc@{}}
\toprule
Operator class name & Tightness \\
\midrule
	\href{https://pepit.readthedocs.io/en/latest/api/operators.html#monotone}{Monotone (maximally)} & \ding{52}\\
	\href{https://pepit.readthedocs.io/en/latest/api/operators.html#strongly-monotone}{Strongly monotone (maximally)} & \ding{52} \\
	\href{https://pepit.readthedocs.io/en/latest/api/operators.html#cocoercive}{Cocoercive} & \ding{52}\\
	\href{https://pepit.readthedocs.io/en/latest/api/operators.html#lipschitz-continuous}{Lipschitz continuous} & \ding{52} \\
    \href{https://pepit.readthedocs.io/en/latest/api/operators.html#negatively-comonotone}{Negative comonotone} & \ding{52}\\
    \href{https://pepit.readthedocs.io/en/latest/api/operators.html#cocoercive-and-strongly-monotone}{Cocoercive and strongly monotone} & \ding{54}\\
    
	\href{https://pepit.readthedocs.io/en/latest/api/operators.html#strongly-monotone-and-lipschitz-continuous}{Lipschitz continuous and strongly monotone} & \ding{54}\\
    \href{https://pepit.readthedocs.io/en/latest/api/operators.html#nonexpansive}{Non-expansive} & \ding{52}\\
    \href{https://pepit.readthedocs.io/en/latest/api/operators.html#linear-operator}{Linear} & \ding{52}\\
    \href{https://pepit.readthedocs.io/en/latest/api/operators.html#symmetric-linear-operator}{Symmetric Linear} & \ding{52}\\
    \href{https://pepit.readthedocs.io/en/latest/api/operators.html#skew-symmetric-linear-operator}{Skew-symmetric Linear} & \ding{52}\\

\bottomrule
\end{tabular}}
\end{center}
\end{table*}
\paragraph{Performance measures and initial conditions.} An important degree of freedom of the package is to allow a large panel of performance measures and initial conditions.
Essentially, everything that can be expressed linearly (or slightly beyond) in function values and quadratically in gradient/iterates (i.e., linear in the Gram representation of Section~\ref{s:example}) might be considered.
Following the notation of Section~\ref{s:code} (and denoting by $x_\star$ an optimal point of $F$), typical examples of initial conditions include:
\begin{itemize}
    \item $\|x_0-x_\star\|^2_2\leqslant 1$,
    \item $\|\nabla F(x_0)\|^2_2\leqslant 1$ (when $F$ is differentiable, otherwise similar criterion involving some subgradient of $F$ might be used),
    \item $F(x_0)-F(x_\star)\leqslant 1$,
    \item any linear combination of the above (see, e.g., examples in the \emph{potential functions} folder of the package).
\end{itemize}
Similarly, typical examples of performance measures (see Table~\ref{tab:perf} for examples) include: $\|x_n-x_\star\|^2_2$, $\|\nabla F(x_n)\|^2_2$ (when $F$ is differentiable), $F(x_n)-F(x_\star)$, or linear combinations and minimum values of the above, \newline
    e.g. $\min_{0\leqslant i\leqslant n} \|\nabla F(x_i)\|^2_2$ (when $F$ is differentiable).
\begin{table*}[ht]
\begin{center}
\caption{Examples of common performance measures. This topic is discussed extensively in excellent references that include~\cite{nemirovsky1992information,nemirovski1994efficient}; see also~\cite[Tables 1--3]{taylor2018exact} for examples in the context of (proximal) gradient descent. Considering appropriate performance measures is key for the analyses, and is particularly exploited when looking for appropriate Lyapunov functions, see, e.g., the related~\cite{lessard2016analysis,taylor19bach,taylor2018lyapunov}.}\label{tab:perf}
{\renewcommand{\arraystretch}{1.8}
\begin{tabular}{@{}lc@{}}
\toprule
Performance measure & Description  \\
\midrule
Distance & $\|x_n-x_\star\|^2_2$  \\
Gradient norm & $\|\nabla F(x_n)\|^2_2$  \\
Function value & $F(x_n)-F(x_\star)$ \\
Contraction & $\|x_n - y_n\|_2^2$  \\
Lyapunov functions & $a(F(x_n)-F(x_\star)) + \begin{pmatrix}x_n-x_\star \\ \nabla F(x_n)\end{pmatrix}^\top (P\otimes I_d) \begin{pmatrix}x_n-x_\star \\ \nabla F(x_n) \end{pmatrix}$ \\
\bottomrule
\end{tabular}
}
\end{center}
\end{table*}

\paragraph{Examples.} \textsc{PEPit} contains about $75$ examples that can readily be used, instantiating the different sets of black-box oracles, problem classes, and initial condition/performance measures. Those examples can be found in the folder \textsc{PEPit/examples/}.

\paragraph{Contributing.} \textsc{PEPit} is designed to allow users to easily contribute to add features to the package. Classes of functions (or operators) as well as black-box oracles can be implemented by following the \href{https://pepit.readthedocs.io/en/latest/contributing.html}{contributing guidelines} from the documentation. We also welcome any new example for analyzing a method/setting that is not already present in the toolbox.

\FloatBarrier

\section{A few additional numerical examples}\label{sec:numerical-examples}

    The following section provides a few additional numerical worst-case analyses obtained through \textsc{PEPit}; namely an accelerated gradient method~\cite{Nest03a}, an accelerated Douglas-Rachford splitting~\cite{patrinos2014douglasrachford}, and point-SAGA~\cite{defazio2016simple} (a proximal method for finite-sum minimization).

    \subsection{Analysis of an accelerated gradient method}\label{subsec:analysis-of-an-accelerated-gradient-method}

        For this example, we focus again on the problem of minimizing a $L$-smooth $\mu$-strongly convex function (problem~\eqref{eq:min_composite} with $F\in\Fmul$).
        We consider a classical accelerated gradient method with constant momentum~\cite{Nest83,Nest03a}.
        It can be described as follows for $t \in \{0, ..., n-1\}$ with $y_0=x_0$:
        \begin{equation}\label{eq:Nest}
            \begin{aligned}
                x_{t+1} & = y_t - \alpha \nabla F(y_t),\\
                y_{t+1} & = x_{t+1} + \beta (x_{t+1} - x_{t}),
            \end{aligned}
        \end{equation}
        with $\kappa = \frac{\mu}{L}$, $\alpha = \frac{1}{L}$ and $\beta = \frac{1 - \sqrt{\kappa}}{1 + \sqrt{\kappa}}$.
        We decide to compute the smallest possible $\tau(n, L , \mu)$ such that the guarantee
        \begin{equation*}
            F(x_n) - F_\star \leqslant \tau(n, L, \mu) \left(F(x_0) -  F(x_\star) + \frac{\mu}{2}\|x_0 - x_\star\|^2_2\right),
        \end{equation*}
        holds for all $d\in\mathbb{N}$, $F\in\Fmul$, $x_0,x_n,x_\star\in\mathbb{R}^d$ where $x_n$ is an output of the accelerated gradient method~\eqref{eq:Nest} and $x_\star$ is the minimizer of $F$.
        In this setting, $\tau(n, L, \mu)$ can be computed as the worst-case value of $F(x_n)-F_\star$ (the performance metric) when $ F(x_0) - F_\star + \frac{\mu}{2}\| x_0 - x_\star\|^2_2\leqslant 1$ (initial condition).
        As a reference, we compare the output of \textsc{PEPit} to the following worst-case guarantee~\cite[Corollary 4.15]{daspremont2021acceleration}:
        \begin{equation}\label{eq:nest_ref}
            F(x_n) - F_\star \leqslant \left(1 - \sqrt{\frac{\mu}{L}}\right)^{n}\left(F(x_0) - F_\star + \frac{\mu}{2}\| x_0 - x_\star\|^2_2\right).
        \end{equation}
        A comparison between the output of \textsc{PEPit} and~\eqref{eq:nest_ref} is presented in Figure~\ref{fig:AGM}, where we see that~\eqref{eq:nest_ref} could be slightly improved to better match the worst-case behavior of the algorithm.
        The corresponding code can be found in \href{https://pepit.readthedocs.io/en/latest/_modules/PEPit/examples/unconstrained_convex_minimization/accelerated_gradient_strongly_convex.html#wc_accelerated_gradient_strongly_convex}{accelerated\_gradient\_strongly\_convex.py} from the directory
        \newline \href{https://pepit.readthedocs.io/en/latest/examples/a.html}{PEPit/examples/unconstrained\_convex\_minimization/}.

    \subsection{Analysis of an accelerated Douglas-Rachford splitting}\label{subsec:analysis-of-an-accelerated-douglas-rachford-splitting}

        In this section, we provide a simple \textsc{PEPit} example for studying an accelerated Douglas-Rachford splitting method.
        This method was introduced in~\cite{patrinos2014douglasrachford} where a worst-case analysis is provided for quadratic minimization.
        We perform a worst-case analysis numerically for a slightly more general setting:
        \begin{equation*}
            F_\star \eqdef \min_x \{F(x) \equiv f_1(x) + f_2(x)\},
        \end{equation*}
        where $f_1$ is closed proper and convex, and $f_2$ is $\mu$-strongly convex and $L$-smooth.
        This section focuses on the following accelerated Douglas-Rachford splitting method, described in~\cite[Section 4]{patrinos2014douglasrachford}:
        \begin{equation*}
            \begin{aligned}
                x_t     &= \mathrm{prox}_{\alpha f_2}(u_t), \\
                y_t     &= \mathrm{prox}_{\alpha f_1}(2x_t - u_t), \\
                w_{t+1} &= u_t +\theta (y_t - x_t), \\
                u_{t+1} &= \left\{\begin{array}{ll} w_{t+1}+\frac{t-1}{t+2}(w_{t+1}-w_t)\, & \mbox{if } t \geqslant 1,\\
                w_{t+1} & \mbox{otherwise,} \end{array} \right.
            \end{aligned}
        \end{equation*}
        where $\mathrm{prox}$ denotes the usual proximal operator, available in \textsc{PEPit} through the operation \textsc{proximal\_step}, as exemplified below.
        Note that we only show the algorithm description here, the full \textsc{PEPit} code for this example can be found in the file \href{https://pepit.readthedocs.io/en/latest/_modules/PEPit/examples/composite_convex_minimization/accelerated_douglas_rachford_splitting.html#wc_accelerated_douglas_rachford_splitting}{accelerated\_douglas\_rachford\_splitting.py} from the directory
        \newline \href{https://pepit.readthedocs.io/en/latest/examples/b.html}{PEPit/examples/composite\_convex\_minimization/}.

\begin{lstlisting}
# Compute n steps of
# An accelerated Douglas-Rachford splitting
for t in range(n):
    x[t], _, _ = proximal_step(u[t], func2, alpha)
    y, _, fy = proximal_step(2*x[t] - u[t],
                             func1, alpha)
    w[t+1] = u[t] + theta * (y-x[t])
    if t >= 1:
        u[t+1] = w[t+1] + (t-1)/(t+2) * (w[t+1]-w[t])
    else:
        u[t+1] = w[t+1]
\end{lstlisting}

        When $f_2$ is a $L$-smooth $\mu$-strongly convex quadratic function, the following worst-case guarantee is provided by~\cite[Theorem 5]{patrinos2014douglasrachford}:
        \begin{equation}\label{eq:ref_quad_ADRS}
            F(y_n) - F_\star \leqslant  \frac{2\|w_0 - w_\star\|^2_2}{\alpha \theta(n+3)^2},
        \end{equation}
        when $\theta = \frac{1-\alpha L}{1 + \alpha L}$ and $\alpha < \frac{1}{L}$.
        A numerical worst-case guarantee for the case $L=1$, $\mu=0.01$, $\alpha=0.9$, $\theta = \frac{1-\alpha L}{1 + \alpha L}$ is provided on Figure~\ref{fig:ADRS} for a few different values of $N$, where we use~\eqref{eq:ref_quad_ADRS} as a reference for comparison.
        For each of those values, \textsc{PEPit} computed a tight (up to numerical precision) worst-case value for which we are not aware of any proven analytical worst-case guarantee beyond the quadratic setting.
        We see that \textsc{PEPit} provides an improvement over this guarantee, even when the problem under consideration is not quadratic.

    \subsection{Analysis of point-SAGA}\label{subsec:analysis-of-point-saga}

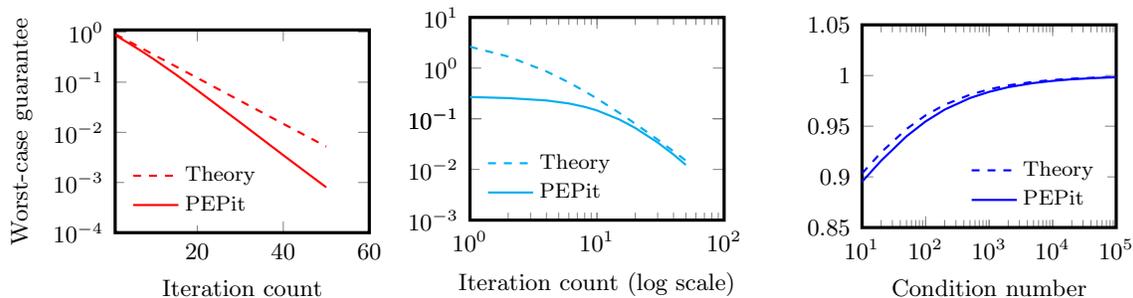
\begin{figure}[!ht]
            \begin{subfigure}[t]{0.31\textwidth}
                \centering
                \begin{tikzpicture}
                    \begin{semilogyaxis}[legend style={draw=none},legend cell align={left}, legend pos={south west},plotOptions4,width=1.05\linewidth,xlabel={Iteration count},xtick={0,20,40,60},xmax=60]
                        \addplot[dashed, red, thick] table [y=theoryfgm, x=condition]{data/fgm_n.txt};
                        \addplot[red, thick] table [y=fgm, x=condition]{data/fgm_n.txt};
                    \end{semilogyaxis}
                \end{tikzpicture}
                \caption{Accelerated gradient method: strong convexity parameter fixed to $\mu=0.1$. Worst-case guarantee on $F(x_n)-F_\star$ as a function of~$n$.}
                \label{fig:AGM}
            \end{subfigure}
            \hfill
            \begin{subfigure}[t]{0.31\textwidth}
                \centering
                \begin{tikzpicture}
                    \begin{loglogaxis}[legend pos=south west,legend style={draw=none},legend cell align={left},plotOptions2,width=1.05\linewidth,ylabel={}]
                        \addplot[dashed, cyan, thick] table [y=theoryadrs, x=condition]{data/adrs_n.txt};
                        \addplot[cyan, thick] table [y=drs, x=condition]{data/adrs_n.txt};
                    \end{loglogaxis}
                \end{tikzpicture}
                \caption{Accelerated Douglas-Rachford splitting with parameter $\alpha = 0.9$ and $\mu=0.1$. Worst-case guarantee on $F(y_n)-F_\star$ as a function of~$n$.}
                \label{fig:ADRS}
            \end{subfigure}
            \hfill
            \begin{subfigure}[t]{0.31\textwidth}
                \centering
                \begin{tikzpicture}
                    \begin{semilogxaxis}[legend pos=south east,legend style={draw=none},legend cell align={left},plotOptions3,width=1.05\linewidth,ylabel={}]
                        \addplot[dashed, blue, thick] table [y=theorypsaga, x=condition]{data/psaga.txt};
                        \addplot[blue, thick] table [y=psaga, x=condition]{data/psaga.txt};
                    \end{semilogxaxis}
                \end{tikzpicture}
                \caption{Point-SAGA: $n=5$ functions. Worst-case guarantee on $E_{j_t}[V(x^{(t+1)}_{j_t}]$ as a function of the condition number $\kappa = \frac{L}{\mu}$.}
                \label{fig:PSAGA}
            \end{subfigure}
            \caption{Comparisons between (numerical) worst-case bounds from \textsc{PEPit}~(\textbf{plain lines}) VS. reference established worst-case guarantees~(\textbf{dashed lines}) for three different optimization methods.
            For simplicity, we fixed smoothness constants to $L=1$.}
            \label{fig:examples}
        \end{figure}

        In this section, we use \textsc{PEPit} for studying point-SAGA~\cite{defazio2016simple}, a stochastic algorithm for finite sum minimization:
        \begin{equation*}
            F_\star \eqdef \min_x \left\{F(x) \equiv \frac{1}{n} \sum_{i=1}^{n} f_i(x)\right\},
        \end{equation*}
        where $f_1, \ldots, f_n$ are $L$-smooth and $\mu$-strongly convex functions with a proximal operator available for each of them.
        At each iteration $t$, point-SAGA picks $j_t\in\{1,\ldots,n\}$ uniformly at random and performs the following updates (a superscript is used for denoting iteration numbers; the subscript is used for referring to the function $f_{j_t}$ chosen uniformly at random):
        \begin{equation*}
            \begin{aligned}
                z_{j_t}^{(t)} &= x^{(t)} + \gamma \left(g_{j_t}^{(t)} - \frac{1}{n} \sum_i g_i^{(t)}\right),\\
                x^{(t+1)}_{j_t} &= \mathrm{prox}_{\gamma f_{j_t}} \left(z_{j_t}^{(t)}\right),\\
                g_{j_t}^{(t+1)} &= \frac{1}{\gamma}\left(z_{j_t}^{(t)} - x^{(t+1)}_{j_t}\right),
            \end{aligned}
        \end{equation*}
        where $\gamma=\frac{\sqrt{(n - 1)^2 + 4n\frac{L}{\mu}}}{2Ln} - \frac{\left(1 - \frac{1}{n}\right)}{2L}$ is the step-size.
        In this example, we use a Lyapunov (or potential / energy) function
        $V(x) = \frac{1}{L \mu}\frac{1}{n} \sum_{i \leqslant n} \|\nabla f_i(x) - \nabla f_i(x_\star)\|^2_2 + \|x - x_\star\|^2_2,$
        and compute the smallest $\tau(n, L, \mu)$ such that the guarantee
        \begin{equation*}
            \mathbb{E}_{j_t}\Big[V\big(x^{(t+1)}_{j_t}\big)\Big] \leqslant \tau(n, L, \mu) V(x^{(t)}),
        \end{equation*}
        holds for all $d\in\mathbb{N}$, $f_i\in\Fmul(\Rd)$ (for all $i=1,\ldots,n$), $x^{(t)}\in\mathbb{R}^d$ where $x^{(t+1)}_{j_t}$ is the (random) output generated by point-SAGA, and the expectation is taken over the randomness of $j_t$ {(and computed in practice as an average of all the possible $j_t$)}.
        The following simple worst-case guarantee is provided in~\cite[Theorem 5]{defazio2016simple} and is used as a reference:
        \begin{equation}\label{eq:ref_ptSAGA}
            \mathbb{E}_{j_t}[V(x^{(t+1)}_{j_t})] \leqslant \frac{1}{1 + \mu\gamma} V(x^{(t)}).
        \end{equation}
        We compare~\eqref{eq:ref_ptSAGA} to \textsc{PEPit}'s tight (up to numerical precision) output in Figure~\ref{fig:PSAGA}.
        We see that the worst-case guarantee~\eqref{eq:ref_ptSAGA} can be slightly improved, although pretty accurate, particularly for large values of the condition number.
        The corresponding \textsc{PEPit} code of this example can be found in the file \href{https://pepit.readthedocs.io/en/latest/_modules/PEPit/examples/stochastic_and_randomized_convex_minimization/point_saga.html#wc_point_saga}{point\_saga.py} from the directory
        \newline \href{https://pepit.readthedocs.io/en/latest/examples/d.html}{PEPit/examples/stochastic\_and\_randomized\_convex\_minimization/}.

        \FloatBarrier

\section{Conclusion}\label{s:ccl}
    The \textsc{PEPit} package, briefly described in this paper, aims at providing simplified access to worst-case analyses of first-order optimization methods in \textsc{python}.
    To achieve this goal, \textsc{PEPit} implements the performance estimation approach while allowing to avoid the possibly heavy semidefinite programming modeling steps.
    The first version of the package already contains about $75$ examples of first-order methods that can be analyzed through this framework.
    Those examples allow either reproducing or tightening, numerically, known worst-case analyses or provide new ones depending on the particular method and problem class at hand.
    
    Overall, we believe that this package allows quick (in)validations of proofs (a step towards reproducible theory) which should help both the development and the review process in optimization.
    We also argue that this is a nice pedagogical tool for learning algorithms together with their worst-case properties just by playing with them.
    Possible extensions under consideration for future versions include an option for searching for Lyapunov (or potential/energy) functions~\cite{taylor19bach,taylor2018lyapunov,upadhyaya2023automated},
    for disproving convergence~\cite{goujaud2023counter}, as well as a numerical proof assistant, and to incorporate recent extensions of PEP/IQCs to distributed and decentralized optimization~\cite{sundararajan2020analysis,colla2021automated}.

\begin{acknowledgements}
    The work of B. Goujaud and A. Dieuleveut is partially supported by ANR-19-CHIA-0002-01/chaire SCAI, and Hi!Paris.
    A. Taylor acknowledges support from the European Research Council (grant SEQUOIA 724063).
    This work was partly funded by the French government under the management of Agence Nationale de la Recherche as part of the ``Investissements d’avenir'' program, reference ANR-19-P3IA-0001 (PRAIRIE 3IA Institute).
    The authors would like to deeply thank the two anonymous referees as well as a technical editor for extremely constructive feedback that contributed to improving this work as well as the related package.
\end{acknowledgements}

\vskip 0.2in

\bibliographystyle{spmpsci}
\bibliography{references}

\end{document}